\documentclass{article}

\usepackage{graphicx} 

\usepackage{framed,multirow}

\usepackage{amsfonts,amssymb,booktabs}
\usepackage{latexsym}
\usepackage{hyperref}

\usepackage{url}
\usepackage{xcolor}
\definecolor{newcolor}{rgb}{.8,.349,.1}
\usepackage{authblk}
\usepackage{amsmath}

\usepackage{xcolor}

\usepackage{algorithmic}
\usepackage{algorithm,amsfonts,amssymb,amsthm}
\usepackage{mathrsfs}
\usepackage{appendix}
\usepackage{ascmac}

\usepackage{cleveref} 
\creflabelformat{equation}{#2(#1)#3}
\crefname{equation}{}{}

\crefname{figure}{Fig.}{}
\crefname{table}{Table}{}
\crefname{section}{Section}{}

\newtheorem{theorem}{Theorem}

\newtheorem{lemma}[theorem]{Lemma}

\newtheorem{remark}[theorem]{Remark}

\def \v{\mathbf{v}}
\def \r{\mathbf{r}}

\def \rp{\mathbf{r}_{p}}
\def \rq{\mathbf{r}_{q}}

\def \n{\mathbf{n}}

\def \nq{\mathbf{n}_q}

\def \v{\mathbf{v}}
\def \r{\mathbf{r}}

\def \rp{\mathbf{r}_{p}}
\def \rq{\mathbf{r}_{q}}

\def \nq{\mathbf{n}_q}
\def \n{\mathbf{n}}

\newcommand{\R}{{\mathbb{R}}}

\title{Efficient Exact Quadrature of Regular Solid Harmonics Times Polynomials Over Simplices in $\R^3$}

\author[1]{Shoken Kaneko\footnote{kaneko60@umd.edu}}
\author[1]{Ramani Duraiswami\footnote{ramanid@umd.edu}}
\affil[1]{Department of Computer Science and Institute for Advanced Computer Studies, University of Maryland, College Park, MD 20742, USA.}

\begin{document}

\maketitle

\begin{abstract}
A generalization of a recently introduced recursive numerical method~\cite{gumerov2023efficient} for the exact evaluation of integrals of regular solid harmonics and their normal derivatives over simplex elements in $\R^3$ is presented.
The original \emph{Quadrature to Expansion} (Q2X) method~\cite{gumerov2023efficient} achieves optimal per-element asymptotic complexity, however, it considered only constant density functions over the elements. 
Here, we generalize this method to support arbitrary degree polynomial density functions, which is achieved in an extended recursive framework while maintaining the optimality of the complexity. The method is derived for 1- and 2- simplex elements in $\R^3$ and can be used for the boundary element method and vortex methods coupled with the fast multipole method.
\end{abstract}

\section{Introduction}

\begin{figure}[htbp]
\centering
\includegraphics[width=7cm]{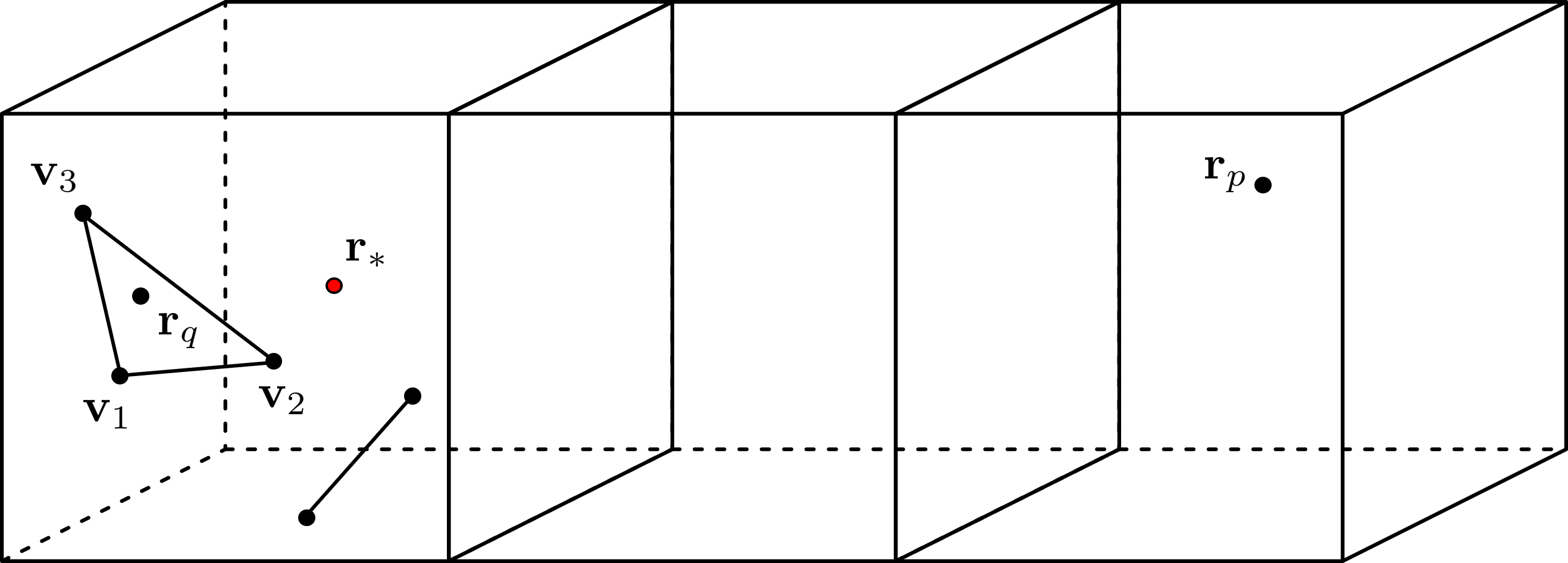}
\vspace{-10.0pt}
\caption{2- and 1-simplex in a cell of an octree data structure commonly used in the FMM-BEM. The evaluation point $\rp$ is located in a well-separated position outside the cell the elements reside.}
\label{fig:octreeCell}
\vspace{-15.0pt}
\end{figure}

Recently, Gumerov et al. introduced \emph{Quadrature to Expansion} (Q2X), a recursive method for the analytical evaluation of integrals of spherical basis functions for the Laplace equation in $\R^3$, aka \emph{solid harmonics}, over $d$-simplex elements for $d\in\{1,2,3\}$ with constant densities~\cite{gumerov2023efficient}. 
This method achieves exact integration of all bases up to truncation degree $p_{\mathrm{s}}$ with optimal complexity $O(p_{\mathrm{s}}^2)$ per element for any $d$ in $\{1,2,3\}$. This is useful for the Boundary Element Method (BEM) \cite{sauter2010boundary} coupled with the Fast Multipole Method (FMM)~\cite{greengard1987fast}. 
The conventional FMM performs an approximate summation of monopoles and dipoles centered at points $\rq$ distributed in space, by expanding them into truncated multipole expansions centered at points $\rp$. Many such expansions are consolidated into one to achieve efficiency. In the Q2X approach, recently introduced by Gumerov et al. \cite{gumerov2023efficient}, surface layer potential integrals over a surface or a line discretized via simplices were represented as such expansions, with the expansion coefficients, obtained by quadrature over the simplex, evaluated analytically via an optimal recursive procedure. 
\cref{fig:octreeCell} shows the geometrical relation of the elements with respect to the expansion center $\r_{*}$ of the solid harmonics, which is typically the centroid of a cell in an octree data structure to which the elements belong.
Recursive analytical methods for high order surface elements have been developed for the close evaluation of layer potential integrals~\cite{kaneko2023recursive}, and such a method is also desired for the integration of regular solid harmonics 
commonly arising in the FMM-BEM. 
\cite{newman1986distributions} discusses a related method for quadrilateral elements, however a full algorithm for computing integrals for all harmonics and all monomial densities with finite degrees has not been described. 
Here, we extend the Q2X method to elements with polynomial densities which allows the evaluation of integrals for all regular solid harmonics up to degree $p_{\mathrm{s}}$ and all monomial density functions up to degree $p_{\mathrm{d}}$ with optimal complexity $O(p_{\mathrm{s}}^2 p_{\mathrm{d}}^d)$ per element. 
All formulation is done for $\R^3$. 

\section{Potential integrals in BEM and vortex methods}
The BEM is widely used for numerical solution of partial differential equations, e.g. the Poisson equation:
\begin{equation}\begin{aligned}
 - \nabla^2 u(\r) = f(\r), \quad \r \in \Omega \subset \R^3,
\end{aligned}\label{eq:HelmholtzEqn} \end{equation}
with field $u$ in domain $\Omega \subset \R^3$, and source $f$ ($=0$ for Laplace). 
The weak form of \cref{eq:HelmholtzEqn} can be written in terms of single- and double layer potentials $L$, $M$~\cite{sauter2010boundary}:
\begin{equation}\begin{aligned}
&\{( c_p \gamma_{0,p} + M \gamma_{0,q} - L \gamma_{1,q}) u\}(\rp) = \{N_0f\}(\rp),\\
&\{L\sigma\}(\rp) \equiv \int_{\rq \in \Gamma} G(\mathbf{r}_p,\mathbf{r}_q) \sigma(\rq) d\Gamma, \\
&\{ M \sigma\}(\rp) \equiv \int_{\rq \in \Gamma} \frac{\partial G(\mathbf{r}_p,\mathbf{r}_q)}{\partial \nq} \sigma(\rq) d\Gamma,\\
\end{aligned}\label{eq:BIE_BF_withRobinBC}
\end{equation}
with $c_p=1/2$ on a smooth boundary, $\n_q$ the outward unit normal, $\gamma_{0}$ and $\gamma_{1}$ the boundary trace and normal derivative operators, respectively, $N_0$ the Newton potential operator, and $G(\rp,\rq)$ the Green function for the Laplace equation:
\begin{equation}\begin{aligned}
&\{ \gamma_{0,q} u \} (\rq) \equiv \lim_{\hat{\r}_q \in \Omega \to \rq \in \Gamma} u(\hat{\r}_q), \quad \rq \in \Gamma=\partial\Omega,\\
&\{ \gamma_{1,q} u \} (\rq) \equiv \nq \cdot \nabla_q u(\rq), \quad \rq \in \Gamma=\partial\Omega, \\
&\{N_0f\}(\rp) = \int_{\rq \in \Omega} G(\rp,\rq) f(\rq) d\Omega, \quad \rp \in \R^3,\\
&G(\rp,\rq) = \frac{1}{4\pi r}, \quad  r \equiv |\mathbf{r}_q-\mathbf{r}_p|.
\end{aligned} \end{equation}

In the BEM the boundary $\Gamma$ is discretized into boundary elements which can be either flat or curved surfaces. The layer potential integrals over these elements are evaluated to form the linear system of equations. 
The densities $\sigma$ are approximated via local, typically polynomial, functions with unknown coefficients which must be determined. 
In the present work we assume the boundary $\Gamma = \partial \Omega$ is a union of boundary elements $\Gamma = \bigcup_i S_i$, where each $S_i$ is a flat triangular element with polynomial density functions. 

Similarly, in vortex methods the Bio-Savart integral is used to compute potentials due to line elements $\Lambda_q$: 
\begin{equation}\begin{aligned}
&\mathbf{H}(\rp) = \int_{\rq \in \Lambda_q}\!\!\!\!\!\!\!\! \frac{\left(\mathbf{I}_q \sigma(\rq)\right) \times\left(\rp-\rq\right)}{4\pi |\rp-\rq| ^3} d \Lambda_q  = \nabla \times (\mathbf{I}_q K_q(\rp))\\
&\{K_q\sigma\}(\rp) \equiv \int_{\rq \in \Lambda_q} G(\mathbf{r}_p,\mathbf{r}_q) \sigma(\rq) d\Lambda_q, \\
\end{aligned}
\end{equation}
where $\mathbf{I}_q$ is the circulation of the $q$-th line element. 

\subsection{Multipole expansion of potentials}
We summarize the formulation also used in \cite{gumerov2023efficient}. The following definition of regular and singular spherical basis functions, $R_{n}^{m}\left( \mathbf{r}\right) $ and $S_{n}^{m}\left( \mathbf{r}\right)$ is accepted:
\begin{equation}\begin{aligned}
&R_{n}^{m}\left( \mathbf{r}\right) = \frac{\left( -1\right) ^{n}i^{\left|
m\right| }}{(n+\left| m\right| )!}r^{n}P_{n}^{\left| m\right| }(\cos \theta
)e^{im\varphi },\\  \label{bas1} 
&S_{n}^{m}\left( \mathbf{r}\right) = i^{-\left| m\right| }(n-\left| m\right|
)!r^{-n-1}P_{n}^{\left| m\right| }(\cos \theta )e^{im\varphi },\quad \quad  \\
&n = 0,1,2,...,\quad m=-n,...,n,\\
&\mathbf{r}  = \left( x,y,z\right)^T =r\left( \sin \theta \cos \varphi
,\sin \theta \sin \varphi ,\cos \theta \right)^T,
\end{aligned}\end{equation}
respectively, where $\left( r,\theta ,\varphi \right) $ are the spherical coordinates of $\mathbf{r}$ and $P_{n}^{m}$ are the associated Legendre
functions~\cite{abramowitz1988handbook}, 
\begin{equation}\begin{aligned}
P_{n}^{m}\left( \mu \right) =\frac{\left( -1\right) ^{m}\left( 1-\mu^{2}\right) ^{m/2}}{2^{n}n!}\frac{d^{m+n}}{d\mu ^{m+n}}\left( \mu^{2}-1\right) ^{n},
\end{aligned}\end{equation}%
for nonnegative $m$. 
$R_n^m$ and $S_n^m$ obey the symmetry:
\begin{equation}
R_{n}^{-m}\left( \mathbf{r}\right) =\left( -1\right) ^{m}\overline{%
R_{n}^{m}\left( \mathbf{r}\right) },\quad S_{n}^{-m}\left( \mathbf{r}%
\right) =\left( -1\right) ^{m}\overline{S_{n}^{m}\left( \mathbf{r}\right) },
\label{bas2.1}
\end{equation}%
where the bar indicates the complex conjugate. 
In these bases the Green's function can be expanded as 
\begin{equation}\begin{aligned}
G\left(\rp, \rq \right) = \sum_{n=0}^{\infty}\frac{(-1)^{n}}{4\pi}\sum_{m=-n}^{n}R_{n}^{-m}\left(\rq-\mathbf{r}_{*}\right) S_{n}^{m}\left( \rp -\mathbf{r}_{*}\right), \label{bas3}
\end{aligned}\end{equation}%
given $\left| \rp -\mathbf{r}_{*}\right| > \left| \rq -\mathbf{r}_{*}\right|$ where $\r_{*}$ is the center of expansion. 
The integrals of the spherical basis functions over element $S_q, \Lambda_q$ are defined as:
\begin{equation}\begin{aligned}
\label{eqn:SLP_and_DLP_mutipoles}
L_{n}^{m}\left( \mathbf{r}_{\ast }\right) &\equiv \frac{(-1)^{n}}{4\pi }\int_{\rq \in S_q}R_{n}^{-m}\left( \mathbf{r}_q -\mathbf{r}_{\ast}\right) \sigma(\rq) dS_q ,\\
M_{n}^{m}\left( \mathbf{r}_{\ast }\right) &\equiv \frac{(-1)^{n}}{4\pi }\mathbf{n}_q\cdot \!\! \int_{\rq \in S_q} \!\!\!\! \!\!\!\! (\nabla R_{n}^{-m}( \mathbf{r}_q -\mathbf{r}_{\ast })) \sigma(\rq) dS_q , \\
K_{n}^{m}\left( \mathbf{r}_{\ast }\right) &\equiv \frac{(-1)^{n}}{4\pi }\int_{\rq \in \Lambda_q}R_{n}^{-m}\left( \mathbf{r}_q -\mathbf{r}_{\ast}\right) \sigma(\rq) d\Lambda_q ,\\
\end{aligned}\end{equation}
and are used to expand the layer potentials for $S_q$ as:
\begin{equation}\begin{aligned}
\{F\sigma\}(\rp) = \sum_{n=0}^{\infty}\sum_{m=-n}^{n} S_n^m(\rp-\r_{*}) F_n^m(\r_{*}),
\end{aligned}\end{equation}
where $F=\{L,M,K\}$.
$R_n^m$ satisfies (see e.g.,~\cite{gumerov2006fast}):
\begin{equation}\begin{aligned}
&\frac{\partial }{\partial \eta }R_{n}^{m}\left( \mathbf{r}\right) =iR_{n-1}^{m+1}\left( \mathbf{r}\right) , \quad    \frac{\partial }{\partial \xi }R_{n}^{m}\left( \mathbf{r}\right) =iR_{n-1}^{m-1}\left( \mathbf{r}\right) ,\\
&\frac{\partial }{\partial z}R_{n}^{m}\left( \mathbf{r}\right) =-R_{n-1}^{m}\left( \mathbf{r}\right) , \!\!\!\! \quad \xi \equiv \frac{x+iy}{2}, \!\!\!\! \quad \eta \equiv \frac{x-iy}{2},\\ 
&\partial_{\eta} \equiv \frac{\partial }{\partial \eta }=\frac{\partial }{\partial x}+i\frac{%
\partial }{\partial y},\quad 
\partial_{\xi} \equiv \frac{\partial }{\partial \xi } = \frac{%
\partial }{\partial x}-i\frac{\partial }{\partial y}.  \label{bas7}
\end{aligned}\end{equation}
Thus, with $\mathbf{n}_q=(n_{x},n_{y},n_{z})^T$ we have:
\begin{equation}\begin{aligned}
\mathbf{n}_q \cdot \nabla R_{n}^{m}\left( \mathbf{r}\right) = i \frac{n_{x}}{2}\left[
R_{n-1}^{m+1}\left( \mathbf{r}\right) + R_{n-1}^{m-1}\left( \mathbf{r}\right) %
\right] \\
+ \frac{n_{y}}{2} \left[ R_{n-1}^{m+1}\left( \mathbf{r}\right)
-R_{n-1}^{m-1}\left( \mathbf{r}\right) \right] -n_{z}R_{n-1}^{m}\left( 
\mathbf{r}\right).  \label{bas8}
\end{aligned}\end{equation}
\section{Problem statement}
For $d=2$, we denote the vertices of the 2-simplex element as $\{\v_1,\v_2,\v_3\}$ and use the parametrization:
\begin{equation}\begin{aligned}
&\mathbf{r}_q\left( u,v\right) = \mathbf{r}_{u}u+\mathbf{r}_{v}v + \mathbf{v}_{1} =A(u,v,0)^T+\v_1= (x,y,z)^T, 
\\ 
&\r_u \equiv \frac{\partial \rq}{\partial u}, \quad \!\! \r_v \equiv \frac{\partial \rq}{\partial v}, \quad \!\!
A \equiv (\r_u,\r_v,\n)^{T}. 
\end{aligned}\end{equation}
In contrast to \cite{gumerov2023efficient} which used constant parametrization, we extend the method to density functions $\sigma$ defined over the elements which are polynomials of the form:
\begin{equation}\begin{aligned}
N(u,v) = \!\!\! \sum_{b=0,c=0}^{b+c\le p_{\mathrm{d}}} A_b^c N_b^c(u,v), \quad \!\! N_b^c(u,v) \equiv u^bv^c.
\end{aligned}\end{equation}
The integrand of the integrals $L_n^m, M_n^m$ have the form:
\begin{equation}\begin{aligned}
Q_{n,b}^{m,c}(u,v) \equiv R_n^m(\rq(u,v)-\r_{*})N_b^c(u,v).
\end{aligned}\end{equation}
We set the origin of the coordinate system so that  $\mathbf{r}_{\ast }=\mathbf{0}$. 
The goal is to derive an algorithm for efficient analytical evaluation of integrals 
\cref{eqn:SLP_and_DLP_mutipoles} for $\sigma=N_b^c(u,v)$:
\begin{equation}\begin{aligned}
\label{eqn:SLP_and_DLP_mutipoles}
L_{n,b}^{m,c}\left( \mathbf{r}_{\ast }\right) &\equiv \frac{(-1)^{n}}{4\pi }\int_{\rq \in S_q}R_{n}^{-m}\left( \mathbf{r}_q -\mathbf{r}_{\ast
}\right) N_b^c(u,v) dS_q ,\\
M_{n,b}^{m,c}\left( \mathbf{r}_{\ast }\right) &\equiv \frac{(-1)^{n}}{4\pi }\mathbf{n}_q\cdot \!\! \int_{\rq \in S_q} \!\!\!\! \!\!\!\! (\nabla R_{n}^{-m}( \mathbf{r}_q -\mathbf{r}_{\ast })) N_b^c(u,v) dS_q , 
\end{aligned}\end{equation}
for all index tuples $(n,m,b,c)$ satisfying $0\le n\le p_{\mathrm{s}}$, $|m| \le n$, $0\le b$, $0\le c$, and $b+c\le p_{\mathrm{d}}$.
The coefficients $K_{n,b}^{m}$ for the case $d=1$ can be expressed by setting $c=0$, $\rq=\r_u u + \v_1$, and replacing $S_q$ by $\Lambda_q$ in $L_{n,b}^{m,c}$.

\section{Q2XP: Q2X for Polynomial elements}
We first describe the case $d=2$. 
The Q2X method~\cite{gumerov2023efficient} was derived by utilizing the fact that the regular spherical basis functions $R_n^m$ are homogeneous polynomials of degree $n$ of arguments $(x,y,z)$. 
Similarly, $N_b^c$ are also homogeneous polynomials of $(u,v)$ with degree $b+c$. 
Hence, by considering $\r$ and $(u,v)$ as functions of $(\xi,\eta,z)$, Euler's homogeneous function theorem gives:
\begin{equation}\begin{aligned}
n R_n^m(\r) & = (\xi \partial_{\xi} + \eta \partial_{\eta} + z \partial_z) R_n^m(\r), \\
(b+c) N_b^c(u,v) &= (\xi \partial_{\xi} + \eta \partial_{\eta} + z \partial_z) N_b^c(u,v).
\end{aligned}\end{equation}
With $a_{s,t}$ the $(s,t)$ entry of $A^{-1}$, $\alpha_s^{(\pm)}\equiv a_{s,1}\pm i a_{s,2}$, 
$\mathbf{a}_s$ the $s$-th row of $A^{-1}$, and 
$\beta_s \equiv \mathbf{a}_s\cdot \v_1$ we obtain:
\begin{equation}\begin{aligned}
&(n+b+c) Q_{n,b}^{m,c} + b \beta_1 Q_{n,b-1}^{m,c} + c \beta_2 Q_{n,b}^{m,c-1}\\
&= \xi (i Q_{n-1,b}^{m-1,c} + b Q_{n,b-1}^{m,c} \alpha_1^{(-)} + c Q_{n,b}^{m,c-1}\alpha_2^{(-)} )  \\
 & + \eta (i Q_{n-1,b}^{m+1,c} + b Q_{n,b-1}^{m,c} \alpha_1^{(+)} + c Q_{n,b}^{m,c-1}\alpha_2^{(+)})  \\
 & + z (- Q_{n-1,b}^{m,c} + b a_{13} Q_{n,b-1}^{m,c} + c a_{23} Q_{n,b}^{m,c-1}).\label{eq:rec_Qnmbc}
\end{aligned}\end{equation}
We define 
$\xi =\xi _{u}u+\xi _{v}v+\xi _{0}$, 
$\eta =\eta _{u}u+\eta _{v}v+\eta _{0}$,
and $z =z_{u}u+z_{v}v+z_{0}$, where the subscripts $u$ and $v$ denote partial derivatives with respect to these variables, and
\begin{equation}\begin{aligned}
&\psi_{n,b}^{m,c} \equiv \!\! \int_{0}^{1} \!\!\! \int_{0}^{1-u} \!\!\!\!\! Q_{n,b}^{m,c}\left( u,v\right) dvdu,\quad \!\!\!
q_{n,b}^{m,c} \equiv Q_{n,b}^{m,c}(1,0),\\
&h_{n,b}^{m,c}\left(u\right) \equiv Q_{n,b}^{m,c}\left( u,1-u\right), \quad \!\!\!\!
j_{n,b}^{m,c} \equiv \!\!  \int_{0}^{1} \!\! h_{n,b}^{m,c}\left( u\right) du.
\end{aligned}\end{equation}

\subsection{Recursions}
\begin{lemma}
The coefficients $\psi_{n,b}^{m,c}$ satisfy the recursion:
\begin{equation}\begin{aligned}
\psi_{n,b}^{m,c} 
=&
\frac{ \xi_0 i \psi_{n-1,b}^{m-1,c}  + \eta_0 i \psi_{n-1,b}^{m+1,c}  - z_0  \psi_{n-1,b}^{m,c} + j_{n,b}^{m,c}}{n+b+c+2}
.\label{eqn:rec_inmbc}
\end{aligned}\end{equation}
\end{lemma}
\begin{proof}
\begin{equation}\begin{aligned}
&u\frac{\partial Q_{n,b}^{m,c}}{\partial u}+v\frac{\partial Q_{n,b}^{m,c}}{\partial v}
= \xi \frac{\partial Q_{n,b}^{m,c}}{\partial \xi }
+\eta \frac{\partial Q_{n,b}^{m,c}}{\partial \eta}
+z\frac{\partial Q_{n,b}^{m,c}}{\partial z} \\
&- \left( \xi _{0}\frac{\partial Q_{n,b}^{m,c}}{\partial \xi }
+\eta _{0}\frac{\partial Q_{n,b}^{m,c}}{\partial \eta }
+z_{0}\frac{\partial Q_{n,b}^{m,c}}{\partial z}\right)  \\
&= (n+b+c) Q_{n,b}^{m,c} + b \beta_1 Q_{n,b-1}^{m,c} + c \beta_2 Q_{n,b}^{m,c-1} \\
&-\xi_0 (i Q_{n-1,b}^{m-1,c} + b Q_{n,b-1}^{m,c} \alpha_1^{(-)} + c Q_{n,b}^{m,c-1}\alpha_2^{(-)} )  \\
 & - \eta_0 (i Q_{n-1,b}^{m+1,c} + b Q_{n,b-1}^{m,c} \alpha_1^{(+)} + c Q_{n,b}^{m,c-1}\alpha_2^{(+)})  \\
 & - z_0 (- Q_{n-1,b}^{m,c} + b a_{13} Q_{n,b-1}^{m,c} + c a_{23} Q_{n,b}^{m,c-1}). \label{eqn:proof_lemma1}
\end{aligned}\end{equation}
By integrating both sides over the triangle $T = \{(u,v)|0\le u, 0\le v, u+v\le1\}$ and applying integration by parts, we obtain \cref{eqn:rec_inmbc}. 
See proof of Lemma 2 in \cite{gumerov2023efficient} for details.
\end{proof}
\begin{lemma}
The coefficients $j_{n,b}^{m,c}$ satisfy the recursion:
\begin{equation}\begin{aligned}
&j_{n,b}^{m,c}=\frac{1}{n+b+c+1} (\\
&(\xi_0+\xi_v) (i j_{n-1,b}^{m-1,c} + b j_{n,b-1}^{m,c} \alpha_1^{(-)} + c j_{n,b}^{m,c-1}\alpha_2^{(-)} )\\ 
&+ (\eta_0+\eta_v) (i j_{n-1,b}^{m+1,c} + b j_{n,b-1}^{m,c} \alpha_1^{(+)} + c j_{n,b}^{m,c-1}\alpha_2^{(+)})\\
&+ (z_0+z_v) (- j_{n-1,b}^{m,c} + b a_{13} j_{n,b-1}^{m,c} + c a_{23} j_{n,b}^{m,c-1}) \\
& - b \beta_1 j_{n,b-1}^{m,c} - c \beta_2 j_{n,b}^{m,c-1} 
+ q_{n,b}^{m,c}).\label{eq:rec_jnmbc}
\end{aligned}\end{equation}
\end{lemma}
\begin{proof}
Replace $Q_{n,b}^{m,c}(u,v)$ in \cref{eqn:proof_lemma1} by $h_{n,b}^{m,c}(u)$, integrate both sides over domain $u \in [0,1]$, and apply integration by parts.
See proof of Lemma 1 in \cite{gumerov2023efficient} for details.
\end{proof}
The recursion for $q_{n,b}^{m,c}$ directly follows from \cref{eq:rec_Qnmbc} by substituting $u=1$ and $v=0$. 
\begin{equation}\begin{aligned}
&q_{n,b}^{m,c} = \frac{1}{n+b+c} ( \\
&(\xi_0+\xi_u) (i q_{n-1,b}^{m-1,c} + b q_{n,b-1}^{m,c} \alpha_1^{(-)} + c q_{n,b}^{m,c-1}\alpha_2^{(-)} )\\ 
&+ (\eta_0+\eta_u) (i q_{n-1,b}^{m+1,c} + b q_{n,b-1}^{m,c} \alpha_1^{(+)} + c q_{n,b}^{m,c-1}\alpha_2^{(+)})\\
&+ (z_0+z_u) (- q_{n-1,b}^{m,c} + b a_{13} q_{n,b-1}^{m,c} + c a_{23} q_{n,b}^{m,c-1})\\
& - b \beta_1 q_{n,b-1}^{m,c} - c \beta_2 q_{n,b}^{m,c-1}).\label{eq:rec_qnmbc}
\end{aligned}\end{equation}
The expansion coefficients $L_{n,b}^{m,c}$ and $M_{n,b}^{m,c}$ then can be computed from the $\psi_{n,b}^{m,c}$ coefficients using $J=|\r_u \times \r_v|$:
\begin{equation}\begin{aligned}
L_{n,b}^{m,c} =& \frac{J}{4\pi }(-1)^{n}\psi_{n,b}^{-m,c}, \quad M_{n,b}^{m,c} =\frac{J}{4\pi }(-1)^{n}l_{n,b}^{-m,c}, \\
l_{n,b}^{-m,c} =& i\frac{n_{x}}{2}\left[ \psi_{n-1,b}^{-m+1,c}+\psi_{n-1,b}^{-m-1,c}\right] \\
&+ \frac{n_{y}}{2}\left[ \psi_{n-1,b}^{-m+1,c}-\psi_{n-1,b}^{-m-1,c}\right] -n_{z}\psi_{n-1,b}^{-m,c}. \label{eq:i_to_L_and_M}
\end{aligned}\end{equation}
\begin{remark}
Q2XP achieves optimal per-element complexity $O(p_{\mathrm{s}}^2 p_{\mathrm{d}}^d)$ for evaluating all index tuples  in question ($(n,m,b,c)$ for $d=2$ and $(n,m,b)$ for $d=1$), whereas a naive application of Gauss-Legendre quadrature requires $O(p_{\mathrm{s}}^2 p_{\mathrm{d}}^d (p_{\mathrm{s}}+p_{\mathrm{d}})^d)$ for exact evaluation. 
\end{remark}

\subsection{Starting values for recursions}
Note that in all recursions $\psi_{n,b}^{m,c}$, $j_{n,b}^{m,c}$,
and $q_{n,b}^{m,c}$ should be set to zero for $\left| m\right| >n$. 
This follows from the fact that $R_{n}^{m}\left( \mathbf{r}\right) = 0$ for $\left| m\right| >n$. 
The starting values are:
\begin{equation}\begin{aligned}
q_{0,b}^{0,c} &= \delta_{c,0},  \quad
j_{0,b}^{0,c} &= \int_0^1  u^b (1-u)^c d u \equiv \kappa_{b,c}\label{eq:starting_qnmbc_jnmbc}
\end{aligned}
\end{equation}
where table $\kappa_{b,c}$ can be computed using recursion:
\begin{equation}
\begin{aligned}
\kappa_{b, c} 
= \frac{c}{b+1} \kappa_{b+1, c-1} \quad  (c>0), \quad
\kappa_{b, 0} =
 \frac{1}{b+1}  \quad (c=0).
\end{aligned}
\end{equation}
Lastly, the starting values for $\psi_{0,b}^{0,c}$ are given by:
\begin{equation}
\begin{aligned}
\psi_{0,b}^{0,c} 
&= \int_0^1 \int_0^{1-u} u^b v^c d v d u 
=\int_0^1u^b \left[\frac{v^{c+1}}{c+1}\right]_0^{1-u}du \\
&= \frac{1}{c+1} \int_0^{1} u^b (1-u)^{c+1} du = \frac{\kappa_{b,c+1}}{c+1} .\label{eq:starting_inmbc}
\end{aligned}
\end{equation}
From \cref{eq:rec_qnmbc}, \cref{eq:rec_jnmbc}, and \cref{eqn:rec_inmbc}
follows \cref{algo_Q2XP}.
\begin{algorithm}[htbp]
\caption{Evaluate multipole expansion coefficients $L_{n,b}^{m,c}$ and $M_{n,b}^{m,c}$ for $0\le n\le p_{\mathrm{s}}$, $|m|\le n$, $0\le b+c\le p_{\mathrm{d}}$}
\begin{algorithmic}\label{algo_Q2XP}
\STATE{1. Compute $q_{0,b}^{0,c}$, $j_{0,b}^{0,c}$, and $\psi_{0,b}^{0,c}$ using \cref{eq:starting_qnmbc_jnmbc,eq:starting_inmbc}}. 
\STATE{2. Compute coefficients $q_{n,b}^{m,c}$ using recursion \cref{eq:rec_qnmbc}}.
\STATE{3. Compute coefficients $j_{n,b}^{m,c}$ using \cref{eq:rec_jnmbc} and $q_{n,b}^{m,c}$}.
\STATE{4. Compute coefficients $\psi_{n,b}^{m,c}$ using \cref{eqn:rec_inmbc} and $j_{n,b}^{m,c}$}.
\STATE{5. Compute coefficients $L_{n,b}^{m,c}$ and $M_{n,b}^{m,c}$ using \cref{eq:i_to_L_and_M}}.
\end{algorithmic}
\end{algorithm}
The case $d=1$ can be easily obtained by minor modifications to the case $d=2$, i.e. setting $c=0$, $\xi_v=\eta_v=z_v=0$, $J=|\r_u|$,  modifying the recursions accordingly, skipping step 4 and 5 in \cref{algo_Q2XP}, and computing the result by: 
\begin{equation}\begin{aligned}
K_{n,b}^{m} =& \frac{J}{4\pi }(-1)^{n}j_{n,b}^{-m,0}. 
\end{aligned}\end{equation}

\section{Numerical experiment}
\begin{figure}[htbp]
\centering
\includegraphics[width=6cm]{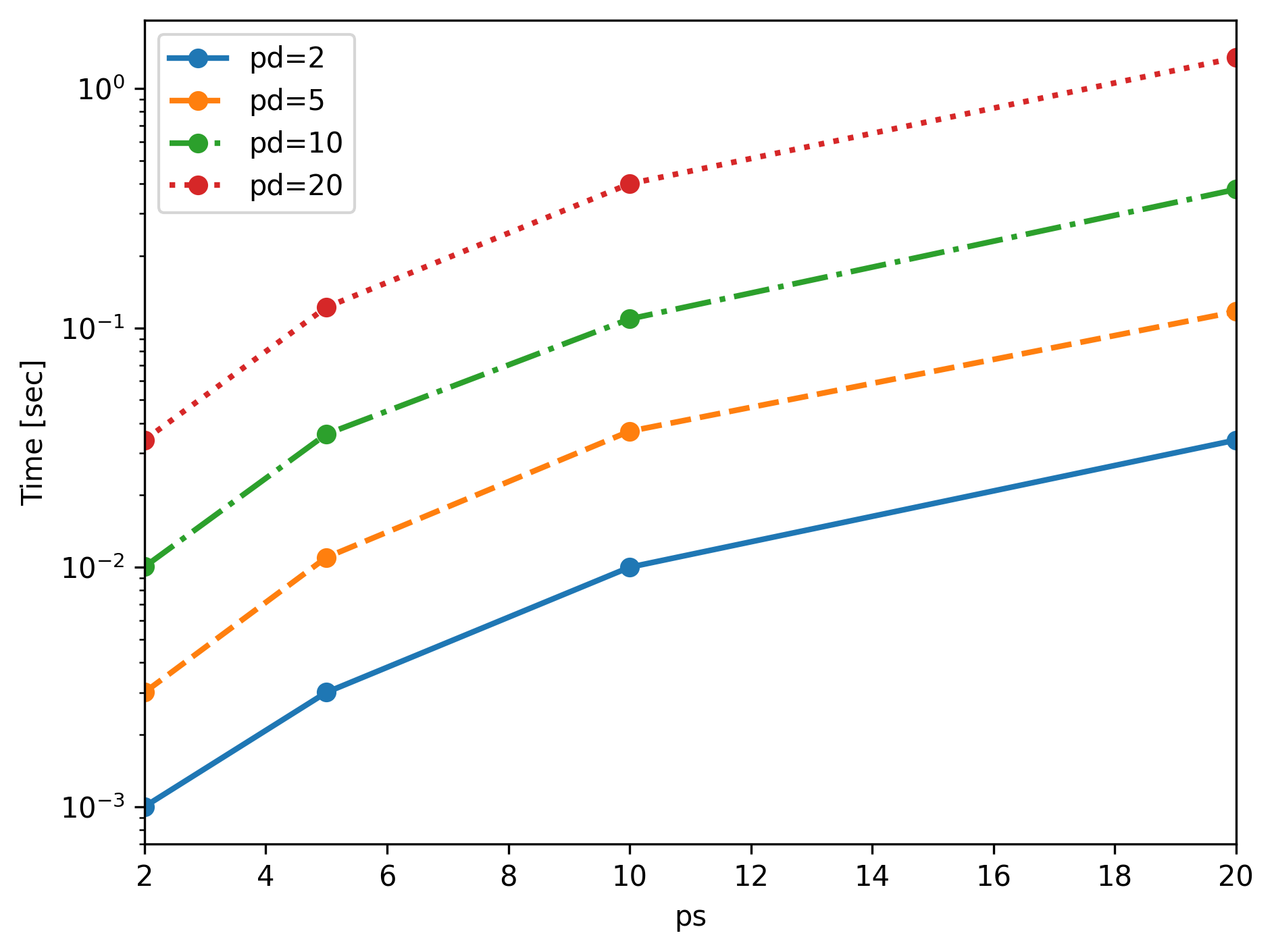}
\vspace{-10.0pt}
\caption{Wall clock times for running Q2XP up to $p_{\mathrm{s}}=p_{\mathrm{d}}=20$.}
\label{fig:times}
\vspace{-10.0pt}
\end{figure}

Q2XP was tested for the case $d=2$ on a configuration given by: 
$\r_c = \left(\sqrt{3}/2,0,0\right)^T$,
$\mathbf{v}_{1} =\mathbf{r}_{c}+r_{t}\left( 1,0,0\right)^T$,
$\mathbf{v}_{2,3}=\mathbf{r}_{c}+r_{t}\left( -1/2,\pm \sqrt{3}/2, 0\right)^T$, and 
$r_t = 0.1$. 
The wall clock times for running a Python implementation of Q2XP for polynomial degrees of up to $p_{\mathrm{s}}=p_{\mathrm{d}}=20$ are shown in \cref{fig:times}.
Given a complexity expression $t=Cp_{\mathrm{s}}^{\alpha}p_{\mathrm{d}}^{\beta}$ with $t$ the computation time and $C$ a constant, the exponents are extracted as $\alpha=1.75$ and $\beta=1.83$ from the numerical experiment. 
The maximum relative difference over the $L_{n,b}^{m,c}$ and $M_{n,b}^{m,c}$ coefficients for truncation degrees $p_{\mathrm{s}}=p_{\mathrm{d}}=10$ computed by exact Gauss-Legendre quadrature and Q2XP was $2.7\times10^{-14}$, indicating good accuracy. 

\section{Conclusion}
We have extended the recursive algorithm presented in \cite{gumerov2023efficient} for the analytical evaluation of integrals of spherical basis functions over high order $d$-simplex elements arising in the FMM-BEM for the Laplace equation in $\R^3$ to the case of polynomial density functions. All multipole expansion coefficients are evaluated analytically up to spherical basis degree of $p_{\mathrm{s}}$ and the density's polynomial degree of $p_{\mathrm{d}}$ with optimal complexity $O(p_{\mathrm{s}}^2p_{\mathrm{d}}^d)$. This complexity was confirmed via numerical experiments. While we limited the discussion to $d\in\{1,2\}$, the case $d=3$ could also be supported by following the formulation presented here and in~\cite{gumerov2023efficient}. 

\section*{Acknowledgments}
This work is supported by Cooperative Research Agreement W911NF2020213 between the University of Maryland and the Army Research Laboratory, with David Hull and Steven Vinci as Technical monitors. 
Shoken Kaneko thanks Japan Student Services Organization and Watanabe Foundation for scholarships.

\bibliographystyle{unsrt}
\bibliography{refs}

\end{document}